\documentclass{amsart}
\usepackage[mathscr]{eucal}
\usepackage[all]{xy}

\theoremstyle{remark}

\title{$K$-equivalence in Birational Geometry}
\author{Chin-Lung Wang}
\thanks{Supported by NSC 89-2115-M-007-045 and NSC 90-2115-M-007-025}
\date{December 18th, 2001. Revised September 30th, 2002.}
\address{Department of Mathematics, National Tsing-Hua University, Hsinchu, Taiwan}
\email{dragon@math.nthu.edu.tw}

\begin{document}
\maketitle

In this article we survey the background and recent development on the
$K$-equivalence relation among birational manifolds. The content is based on
the author's talk at ICCM-2001 at Taipei. I would like to dedicate this article
to Professor Chern, Shiing-shen to celebrate his 90th birthday. For manifolds,
$K$-equivalence is the same as $c_1$-equivalence. In this sense, a major part
of birational geometry is really to understand the geometry of the first Chern
class.\par

After a brief historical sketch of birational geometry in \S1, we define in \S2
the $K$-partial ordering and $K$-equivalence in a birational class and discuss
geometric situations that will lead to these notions. One application to the
filling-in problem for threefolds is given. In \S3 we discuss motivic aspect of
$K$-equivalence relation. We believe that $K$-equivalent manifolds have the
same Chow motive though we are unable to prove it at this moment. Instead we
discuss various approaches toward the corresponding statements in different
cohomological realizations. \S4 is devoted to the {\it Main Conjectures} and
the proof of a weak version of it. Namely, up to complex cobordism,
$K$-equivalence can be decomposed into composite of classical flops. Finally in
\S5 we review some other current researches that are related to the study of
$K$-equivalence relation.\par

\section{A Brief History of Birational Geometry}

If not specifically stated, the ground field is assumed to be the complex
numbers $\mathbb{C}$. Two algebraic varieties are called birational if they
have an isomorphic Zariski open subset. This is equivalent to say that they
have isomorphic rational function fields over the ground field. One of the main
goals in birational geometry is to find a good geometric model that is
convenient for the study of the given algebraic variety or its function
field.\par

\subsection{Minimal Models for Surfaces} (c.f.\ \cite{BPV})
Already at the beginning of the 20th century, Italian algebraic geometers had
adapted the above point of view and successfully applied it to the
classification theory of algebraic surfaces. This led to the famous {\it
Enrique classification}. It started with the {\it Castelnuovo's contraction
theorem}: if a smooth surface $X$ contains at least one $(-1)$ rational curve
$C$, that is $C\cong\mathbb{P}^1$ with $C^2=-1$, then $C$ can be contracted to
a smooth point under a projective birational morphism $\phi_C: X\to X'$ and one
obtains a simplified smooth surface $X'$. By repeating this process in finite
steps, one ends up with a smooth surface called a {\it minimal model}.\par

When $\kappa(X)=-\infty$, that is $\Gamma(X, K_X^m)= 0$ for all
$m\in\mathbb{N}$, {\it Enrique's theorem} says that $X$ is birational to a
ruled surface $C\times\mathbb{P}^1$.\par

When $\kappa(X)\ge 0$, that is $\Gamma(X, K_X^m)\ne 0$ for some
$m\in\mathbb{N}$, the minimal model is unique and it admits semi-ample
canonical bundle. Namely, the minimality of $X$ implies the {\it abundance
theorem}: the pluri-canonical system $\Gamma(X,K_X^m):X\to \mathbb{P}^N$ is a
morphism for $m$ large. Based on this, a detailed classification of algebraic
surfaces according to the Kodaira dimension is then achieved. \par

Later on, Enrique's classification was extended to compact complex surfaces by
Kodaira and to algebraic surfaces over fields of characteristic $p>0$ by
Bombierie and Mumford.\par

\subsection{Minimal Models for Threefolds} (c.f.\ \cite{KMM} \cite{KoMo2})
In 1982, Mori proved the three dimensional generalization of Castelnuovo's
contraction theorem \cite {Mori1}. Mori's theorem opened the way to continue
the {\it minimal model program} (MMP for short). He showed that if a three
dimensional projective manifold $X$ admits a curve $\tilde{C}$ such that
$K_X.\tilde{C}<0$ ($K_X$ is not nef), then by degenerating $\tilde{C}$ inside
$X$ if necessary, one is able to find a rational curve $C$ and a projective
morphism $\phi_R:X \to X'$ associated to the {\it extremal ray}
$R=\mathbb{R}[C]$ such that a curve $C'$ is contracted to a point by $\phi_R$
if and only if $[C']\in R$. He also classified all the possible types of
contractions $\phi_R: X \to X'$. The situation that differs from the surface
case is, the resulting threefold $X'$ is usually mildly singular, with the so
called {\it terminal singularities}.
\par

A priori this seems to be an obstacle to continue the process. Fortunately the
corresponding contraction theorem for terminal threefolds, as well as its
higher dimensional generalization, was soon proved by Kawamata and Shokurov.
However, worse singularities may sometimes occur so that $K_{X'}$ is not even
rationally defined as a line bundle ($X$ is not $\mathbb{Q}$-Gorenstein). This
is the case precisely when the contraction $\phi_R$ is {\it small}, that is,
the exceptional loci of $\phi_R$ has codimension at least two in $X$. In
dimension three, this means that the extremal curve $C$ is isolated. The most
striking idea, which is also the most difficult step in the three dimensional
MMP, is to develop a special type of algebraic surgeries called {\it flips}. A
flip will correct such a pair $C\subset X$ in codimension two along $C$ into
another pair $C' \subset X'$ so that $K_{X'}.C' > 0$, hence it will avoid small
extremal contractions that leads to uncontrollable singularities. \par

The existence of flips in dimension three was proved by Mori in 1988
\cite{Mori3} and the three dimensional MMP was thus established. It states
that, starting with a smooth threefold $X$, after a finite number of divisorial
extremal contractions and/or flips, one ends up with a $\mathbb{Q}$-factorial
terminal model $X'$ which is either (in case $\kappa(X)=-\infty$) a Mori fiber
space $\phi: X' \to S$ under a further extremal contraction of fiber type (with
each fiber $X'_s$ a $\mathbb{Q}$-Fano variety of Picard number $1$) or (in case
$\kappa(X)\ge 0$) a minimal model in Mori's sense, namely $K_{X'}$ is nef. In
the later case, the abundance conjecture that $K_{X'}$ is indeed semi-ample was
subsequently proved by Miyaoka and Kawamata \cite{Kawa2}. Detailed
classifications of threefolds based on the MMP is now of current research
interest.\par

\subsection{Birational Minimal Threefolds and Flops}
In higher dimensions, there are some recent approaches toward the existence of
flips (e.g.\ Shokurov's work for fourfold log-flips), but these have not yet
been completely justified at the time of this writing. Besides the existence
problem, even in the three dimensional case, minimal models are in general not
unique. Based on Reid and Mori's classification theory of three dimensional
singularities \cite{Mori2}, Koll\'ar and Mori in 1989 \cite{Kollar} and then in
1992 \cite{KoMo1} completely understood the precise relation between two
($\mathbb{Q}$-factorial terminal) birational minimal models. The birational map
can always be decomposed into a finite sequence of algebraic surgeries called
{\it flops}. Roughly speaking, each flop is obtained by removing one chain of
rational curves $C$ (corresponding to certain Dynkin diagram) in $X$ with
$K_X|_C=0$ then gluing back $C$ into the open space $X \backslash C$ in a
different manner. \par

This statement was first shown by Kawamata in 1986 \cite{Kawa1}. Koll\'ar and
Mori's method has the advantage to complete the classification of three
dimensional flops (and also flips in 1992), and hence showed that although
three dimensional birational minimal models are in general not homotopically
equivalent, they do have naturally isomorphic ordinary (resp.\ intersection)
cohomology groups, mixed (resp.\ intersection pure) Hodge structures, set of
germs of isolated singularities and local moduli spaces (they actually showed
that flops can be performed simultaneously in flat families). These results
have put the three dimensional minimal model theory into a solid and useful
stage.\par

\section{$K$-partial Ordering in a Birational Class}

\subsection{$K$-partial Ordering} (c.f.\ \cite{Wang2})
We start with a simple observation of the MMP from the point of view of
canonical divisors. For a birational map $f:X \dashrightarrow X'$ between two
$\mathbb{Q}$-Gorenstein varieties, we say that $X\le_K X'$ (resp.\ $X <_K X'$)
if there is a birational correspondence $(\phi,\phi'): X\leftarrow Y \to X'$
extending $f$ with $Y$ smooth, such that $\phi^*K_X \le_{\mathbb{Q}}
{\phi'}^*K_{X'}$ (resp.\ $<_{\mathbb{Q}}$) as divisors. These relations are
easily seen to be independent of the choice of $Y$. Notice that $X\le_K X'$ and
$X\ge_K X'$ imply $X=_K X'$, that is $\phi^*K_X =_{\mathbb{Q}}
{\phi'}^*K_{X'}$. In this case, we say that $X$ and $X'$ are $K$-equivalent. In
this $K$-partial ordering, divisorial contractions and flips will decrease its
$K$-level while flops inducing $K$-equivalence. It is easy to see that
$K$-equivalent terminal varieties are isomorphic in codimension one. In fact,
more is true in general (Theorem 1.4 in \cite{Wang2}):\par
\smallskip
{\it Let $f: X \dashrightarrow X'$ be a birational map between two varieties
with canonical singularities. Suppose that the exceptional locus $Z\subset X$
is proper and that $K_X$ is nef along $Z$, then $X\le_K X'$. Moreover, If $X'$
is terminal then ${\rm codim}_X Z\ge 2$.}\par
\smallskip
Let us recall the proof briefly. Let $(\phi,\phi'):Y\to X\times X'$
be a resolution of $f$ so that the union of the exceptional set of $\phi$ and
$\phi'$ is a normal crossing divisor of $Y$. Let
$K_Y=_{\mathbb{Q}} \phi^*K_X + E =_{\mathbb{Q}} {\phi'}^*K_{X'} + E'$. So
$$
{\phi'}^*K_{X'} =_{\mathbb{Q}} \phi^* K_X + F, \quad\mbox{with $F:=E - E'$}.
$$
It suffices to show that $F\ge 0$. Let $F = \sum_{j=0}^{n-1} F_j$ with $\dim
\phi'({\rm Supp}\,F_j)=j$. We will show that $F_j\ge 0$ for
$j=n-1,n-2,\cdots,1,0$ inductively. As $E'$ is $\phi'$-exceptional, $F_{n-1}\ge
0$ is clear. Suppose that we have already shown that $F_j\ge 0$ for $j\ge k +
1$.\par

Consider the surface $S_k := H^{n-2-k}.{\phi'}^*L^k$ on $Y$ where $H$ is very
ample on $Y$ and $L$ is very ample on $X'$. We get a relations of divisors on
$S_k$:
$$
{\phi'}^*K_{X'}|_{S_k}=_{\mathbb{Q}}\phi^*K_X|_{S_k} + a-b,
$$
where $H^{n-2-k}.{\phi'}^*L^k.F=a-b$ with both $a$ and $b$ effective. Notice
that $b$ can only come from $F_k$ since $\sum_{j\ge k+1}F_j\ge 0$ and $L^k\cap
\phi'(F_j)=\emptyset$ for $j<k$. Now we look at
$$
b.{\phi'}^*K_{X'}=_{\mathbb{Q}} b.\phi^*K_X + b.a -b^2.
$$
The left hand side is always zero since $\phi'(b)\subset L^k\cap\phi'(F_k)$ is
zero dimensional. Moreover, since ${\phi'}^*K_{X'} =_{\mathbb{Q}} \phi^* K_X$
on $\phi^{-1}(X\backslash Z)$, we must have that $\phi({\rm Supp}\, F)\subset
Z$. In particular, $b.\phi^*K_X \ge 0$. It is also clear that $b.a \ge 0$.
However, if $b \ne 0$ then it is a nontrivial combination of $\phi'$
exceptional curves in $S_k$. By the Hodge index theorem for surfaces we then
have that $b^2<0$, a contradiction. So $b=0$ and $F_k\ge 0$. The codimension
statement is easy and we omit its proof.\par

As a corollary, birational minimal models, if they exist, are $K$-equivalent
and reach the lowest $K$-level among terminal (or even canonical) varieties
within their birational class. It is suggestive to make use of this
$K$-equivalence (quasi-uniqueness) and the minimum property to study minimal
models. \par

\subsection{Filling-in Problem in Dimension Three} (c.f.\ \cite{Wang1})
Among applications of the {\it Mori theory}, we mention only one example which
also makes use of $K$-equivalence relation for fourfolds and an extra technique
called {\it symplectic deformations} that will be important in 4.3 in
formulating the Main Conjectures.\par
\smallskip
{\it Let $\mathscr{X}\to\Delta$ be a projective smoothing
of a nontrivial Gorenstein minimal threefold $\mathscr{X}_0$ over the unit disk. Then,
up to any finite base change,
$\mathscr{X}\to\Delta$ is not $\Delta$-birational to a projective smooth family
$\mathscr{X}'\to\Delta$ of minimal threefolds.}\par

\begin{proof}
By an application of the Shokurov-Koll\'ar connectedness theorem one may show
that $\mathscr{X}$ has at most terminal singularities. Then by 2.1, any
$\Delta$-birational map $f$ will induce $K$-equivalence of $\mathscr{X}$ and
$\mathscr{X}'$. Hence they are isomorphic in codimension one and $f$ induces a
birational map $f_0:\mathscr{X}_0\dashrightarrow \mathscr{X}'_0$ between
minimal models. By Koll\'ar's result on birational minimal threefolds in 1.3,
$\mathscr{X}_0$ can not be $\mathbb{Q}$-factorial since $\mathscr{X}'_0$ is
smooth. Now by a result of Kawamata \cite{Kawa1} (or by the three dimensional
MMP), there is a $\mathbb{Q}$-factorialization $\phi:Y\to \mathscr{X}_0$ with
$\phi$ a small morphism. Again this implies that $Y$ is smooth and there is a
birational map of smooth minimal threefolds $f_0\circ\phi:Y\dashrightarrow
\mathscr{X}'_0$. In particular, $$ H^k(Y)\cong H^k(\mathscr{X}'_0)\cong
H^k(\mathscr{X}'_t)\cong H^k(\mathscr{X}_t) \quad\mbox{for all $k\ge 0$ and
$t\ne 0$}. $$ Now we are in a {\it small contraction/smoothing} diagram: $$
\xymatrix{ Y\ar[d]_\phi & &\\ \mathscr{X}_0\,\ar@{^{(}->}[r]
&\mathscr{X}\ar@{<-^{)}}[r] &\,\,\mathscr{X}_t} $$ In case that $\mathscr{X}_0$
has only ODP's, a simple Mayer-Vietoris argument shows that this is impossible.
In fact, consider a diagram as above in the $C^\infty$ category such that near
each singular point of $\mathscr{X}_0$ it is a small contraction/smoothing
diagram of a germ of ODP. Let $C_i$'s be the rational curves contracted to
those ODP's and let $e:\bigoplus_i \mathbb{Z}[C_i]\to H_2(Y,\mathbb{Z})$ be the
class map, then $H_2(\mathscr{X}_t)={\rm coker}\, e$. So,
$H_2(\mathscr{X}_t)\cong H_2(Y)\Rightarrow {\rm Image}\, e = 0$, which is
impossible since $Y$ is projective.\par

In general, by Reid's classification \cite{Reid}, three dimensional Gorenstein
terminal singularities are exactly isolated cDV points (one parameter
deformations of surface ADE singularities). By Friedman's result
\cite{Freidman}, if $p\in V$ is a germ of an isolated cDV point and $C\subset
U$ is the corresponding germ of the (possibly reducible) exceptional curve
contracted to $p$, then the versal deformation spaces ${\rm Def}(p,V)$ and
${\rm Def}(C,U)$ are both smooth and there is an inclusion map of complex
spaces ${\rm Def}(C,U)\hookrightarrow {\rm Def}(p,V)$. Moreover, one can deform
the complex structure of a small neighborhood of $C$ so that in this new
complex structure $C$ deforms into several $\mathbb{P}^1$'s and the contraction
map deforms to a nontrivial contraction of these $\mathbb{P}^1$'s down to
ODP's, while keeping a neighborhood of these ODP's to remain in the versal
deformations of the germ $p\in V$. \par

We can preform this analytic process for all $C$'s and $p$'s simultaneously in
each corresponding small neighborhoods and then patch them together smoothly.
(In fact, one may keep the overlapped region admitting nearby almost complex
structure which is tamed by the original symplectic form \cite{Wilson}. We call
this {\it locally holomorphic} symplectic deformations.) As a result, we obtain
a deformed $C^\infty$ diagram which satisfies the conditions stated above,
which again leads to a contradiction to the equality of $H^2$.
\end{proof}

This allows us to construct counterexamples to the so called three dimensional
{\it filling-in problem}. In fact, Clemens had constructed $A_2$ degenerations
of (simply connected) quintic Calabi-Yau hypersurfaces in $\mathbb{P}^4$ over
$\Delta$ such that the algebraic sub-family over the punctured disk is
$C^\infty$ trivial (so that one may replace the special fiber, which is a
singular Calabi-Yau with an $A_2$ singularity, by a real six dimensional smooth
manifold to obtain a smooth family over $\Delta$). However, we just show that
this smooth replacement can not be achieved in the algebraic category. \par

One of the key points in the above proof is that birational smooth minimal
threefolds have the same Betti numbers. This motivated the author to consider
the validity of {\it equivalence of Betti numbers} in higher dimensions. It is
clear that one should find methods independent of the MMP to study relations
between $K$-equivalent varieties or manifolds.\par

\subsection{Integration Formalism/A Meta Theorem} (c.f.\ \cite{Wang2})
Starting with a birational correspondence with smooth $Y$:
$$
\xymatrix{&Y\ar[ld]_\phi\ar[rd]^{\phi'}\\X&&X'}
$$
From $K_Y=\phi^*K_X + E$ and $K_Y={\phi'}^*K_{X'} + E'$, we see that $X=_K X'$
is the same as saying that $\phi$ and $\phi'$ have the same {\it holomorphic
Jacobian factor} $E=E'$. If for some geometric/topological invariant $I(X)$
that can be computed from certain integration theory whose {\it change of
variable formula} respects the holomorphic Jacobian factor, then one may
conclude that $X=_K X' \implies I(X)=I(X')$ via
$$
I(X)=\int_X d\mu_X=\int_Y J(E)d\mu_Y = \int_{X'}d\mu_{X'}=I(X').
$$
We are going to discuss several examples of this {\it Meta Theorem}
in the following sections.\par

\section{$K$-equivalence Relation and Motives}

Grothendieck's theory of {\it motives} is in principle {\it the} universal
cohomology theory which admits various realizations as usual cohomologies
(Betti, de Rham, Hodge, $\ell$-adic \'etale and others). The category of
motives is supposed to be a homomorphic image of the category of varieties and
should have many expected linear structures (like Hodge filtration and Galois
actions). Unfortunately such a category has not yet been constructed. The
closest one seems to be the Chow motives, or classical motives. Roughly, this
category has all varieties as its objects and the morphisms ${\rm Hom}_{\rm
motive}(X,X')$ are given by correspondences (cycles) $\Gamma\in A^*(X\times
X')$ modulo an {\it adequate equivalence relation} (e.g.\ rational equivalence,
homological equivalence or numerical equivalence). See e.g.\ \cite{Fulton} for
some basic properties.\par

We would like to convince the reader that for $K$-equivalent manifolds under
birational map $f:X\dashrightarrow X'$, there is a {\it naturally attached}
correspondence $T\in A^{\dim X}(X\times X')$ of the form $T = \bar\Gamma_f +
\sum_i T_i$ with $\bar\Gamma_f\subset X\times X'$ the cycle of graph closure of
$f$ and with $T_i$'s being certain degenerate correspondences (i.e.\ $T_i$ has
positive dimensional fibers when projecting to $X$ or $X'$) such that $T$ is an
isomorphism of Chow motives. Currently we do not know how to prove it but some
statements in various realizations do admit proofs along the line of our
integration formalism.\par

\subsection{Classical Integration and $L^2$ Cohomology}
The first clue for the author to believe in a close relationship between
$K$-equivalent manifolds is a somewhat n\"aive yet exciting idea which involves
degenerate K\"ahler metrics and $L^2$ cohomology. \par

Let $X$ and $X$ be smooth projective (be K\"ahler is enough) and let
$(\phi,\phi'): Y\to X\times X'$ be the birational correspondence which leads to
$X =_K X'$. We may select arbitrary K\"ahler metrics $\omega$ and $\omega'$
with volume 1 on $X$ and $X'$ respectively. Then we pull backs them to $Y$ to
get two degenerate K\"ahler metrics $\phi^*\omega$ and ${\phi'}^*\omega'$. From
$c_1$-equivalence we see that (let $\dim X = n$)
$$
(-\phi^*\partial\bar\partial\log\omega^n) -
(-{\phi'}^*\partial\bar\partial\log{\omega'}^n)
=\partial\bar\partial f
$$
for some $C^\infty$ function $f$ up to a constant. This simplifies to
$({\phi'}^*\omega')^n=e^f(\phi^*\omega)^n$. That is, the two degenerate metrics
$\phi^*\omega$ and ${\phi'}^*\omega'$ have quasi-equivalent volume forms (both
volume forms have the same rate of degeneracy along the common degenerate loci
$E\subset Y$).\par

By using cohomology of $L_2$ smooth differential forms with respect to a
possibly degenerate smooth K\"ahler metric, $H^k(X)\cong
L_2^k(X,\omega)=L_2^k(Y,\phi^*\omega)$. If we may {\it rotate} $\phi^*\omega$
to ${\phi'}^*\omega'$ through (not necessarily K\"ahler) degenerate metrics
$g_t$, $t\in [0,1]$ while keeping the volume degeneracy unchanged, then the
theory of $L^2$ cohomology will lead to a proof of the equivalence of
cohomology groups. \par

One candidate for this rotation is to solve a family of complex Monge-Amper\`e
equations via Yau's solution to the Calabi conjecture \cite{Yau}:
$$
(\tilde\omega + \partial\bar\partial\varphi_t)^n =
e^{t(f+c(t))}(\phi^*\omega)^n,
$$
where $\tilde\omega$ is an arbitrary K\"ahler metric with volume 1 on $Y$ and
$c(t)$ is a normalizing constant at time $t$ to make the right hand side has
total integral 1 over $Y$. At this moment there are still analytical
difficulties of this differential geometric approach that need to be overcome.
\par

\subsection{$p$-adic Integration and \'Etale Cohomology}
In \cite{Wang2}, the author applied the idea of quasi-equivalent volume
elements in the theory of $p$-adic integrals. This extended an earlier result
of Batyrev \cite{Batyrev1} on the equivalence of Betti numbers for birational
Calabi-Yau manifolds to general $K$-equivalent manifolds. In particular this
applies to birational smooth minimal models (c.f.\ \S2.1). \par

We will assume that $X$ and $X'$ are smooth projective. Take an integral model
of the $K$-equivalence diagram, e.g.\ $X\to {\rm Spec}\,S$ etc.\ with $S$ a
finitely generated $\mathbb{Z}$-algebra. For {\it almost all} maximal ideals
$P$ in $S$, in fact Zariski open dense in the maximal spectrum of $S$, we have
good reductions of $X$, $X'$, $Y$, $\phi$ and $\phi'$. In such cases, let
$R=\hat{S}_{P}$, the completion of $S$ at $P$ with residue field
$R/P\cong\mathbb{F}_q$, $q=p^r$ for some $r$. For ease of notation, we use the
same symbol to denote the corresponding object over ${\rm Spec}\,R$. Let
$U_i$'s be a Zariski open cover of $X$ such that $K_X|_{U_i}$ is trivial for
each $i$. Then for a compact open subset $S\subset U_i(R)\subset X(R)$, we
define its $p$-adic measure by
$$
m_X(S)\equiv \int\nolimits_S |\Omega_i|_p.
$$
This is independent of the choice of $\Omega$. The $p$-adic measure of $X(R)$
and $X'(R)$ are the same by the {\it change of variable formula} and $X=_K X'$.
By a direct extension of Weil's formula \cite{Weil}, we see that (let $\dim X =
n$ and $\bar X$ be the special fiber)
$$
m_X(X(R)) = \frac{|\bar X(\mathbb{F}_q)|}{q^n}.
$$
By applying this to finite extensions of $\mathbb{F}_q$, we conclude that $X$
and $X'$ have the same local zeta functions for almost all maximal ideals
$P$.\par

Knowing this for one $P$ already allows us to apply Grothendieck-Deligne's
solution to the celebrated Weil conjecture \cite{Deligne1} to conclude that
$K$-equivalent manifolds have the same Betti numbers. \par

In fact more is true \cite{Wang4}. For simplicity let us assume that $X$, $X'$
and $X=_K X'$ are defined over a $\mathbb{Z}$-algebra $S$ such that the
quotient field $F$ of $S$ is a number field. The general case can be reduced to
the number field case by standard tricks, e.g.\ by taking an $F$-valued point
$\eta:{\rm Spec}\,F\to {\rm Spec}\,S$ and considering the fiber diagram over
$\eta$. We then know that $X$ and $X'$ have the same local zeta functions for
almost all, hence all but finite, $P\in {\rm Spec}\, S$. \par

By the C\v{e}botarev density theorem \cite{Serre}, this implies that, for
suitable prime $p$, the two rational Galois representations $H^k_{et}(X_{\bar
F},\mathbb{Q}_p)$ and $H^k_{et}(X'_{\bar F},\mathbb{Q}_p)$ have isomorphic
semi-simplifications as ${\rm Gal}(\bar F/F)$ modules. In other words, $X$ and
$X'$ have the same ``motives'' in the sense of $L$ functions.
\par

By the {\it Hodge-Tate decomposition theorem} proved by Fontaine and Messing
\cite{FoMe} under certain restrictions or by Faltings' complete version of {\it
$p$-adic Hodge theory} \cite{Faltings}, this then implies the equivalence of
$\mathbb{Q}_p$ (and hence $\mathbb{Q}$) Hodge numbers.\par

More precisely, select a prime $P$ so that $X$ and $X'$ have good reductions
(this is in fact unnecessary) and let $K = F_P$ be the completion with residue
field $k$ of characteristic $p$. Let $G={\rm Gal}(\bar K/K)$ and $\mathbb{C}_p$
be the completion of $\bar K$. Then there exists a natural $G$-equivariant
isomorphism \cite{Faltings}:
$$
\bigoplus\nolimits_i \big(\mathbb{C}_p\otimes_K H^{m-i}(X_K,\Omega^i)(-i)\big)\cong
\mathbb{C}_p\otimes_{\mathbb{Q}_p} H^m_{et}(X_{\bar K},\mathbb{Q}_p),
$$
where $G$ acts on $H^{m-i}(X_K,\Omega^i)$ trivially and on the right hand side
diagonally, and $(j)$ is the Tate twist by $j$-th power of cyclotomic
character. Since $\mathbb{C}_p^G=K$ and $\mathbb{C}_p(i)^G=0$ for $i\ne 0$,
elementary manipulation shows that
$$
h^{i,m-i}= \dim_K \big(\mathbb{C}_p\otimes_{\mathbb{Q}_p}
          H^m_{et}(X_{\bar K},\mathbb{Q}_p)^{ss}(i) \big)^G.
$$\par

Finally we plug in $H^m_{et}(X_{\bar K},\mathbb{Q}_p)^{ss}\cong
H^m_{et}(X'_{\bar K},\mathbb{Q}_p)^{ss}$, which holds by base change theorem,
to conclude the non-canonical equivalence of $\mathbb{Q}_p$-Hodge structures.
\footnote{I am grateful to C.-L.\ Chai and J.-K.\ Yu for discussions on the
sufficiency for determining Hodge numbers from semi-simplifications.
\par In the Algebraic Geometry Conference for Iitaka's 60 at Tokyo, February 2002,
T.\ Ito informed the author that he also obtained the same proof independently
\cite{Ito}.}

\subsection{Grothendieck Group of Varieties and the Hodge Realization}
Let $K_0({\rm Var}_{\mathbb{C}})$ be the Grothendieck ring of complex varieties
with reduced structures. That is, we modulo out the {\it motivic relation}
$[X\backslash Z]=[X]-[Z]$ whenever $Z\subset X$ is a closed subvariety. (For
$X$ a variety, we use $[X]$ to denote its class in $K_0({\rm
Var}_{\mathbb{C}})$.) To see that $K_0({\rm Var}_{\mathbb{C}})$ contains enough
information, notice that the Hodge realization functor $h^{p,q}$ defined on
smooth projective varieties has an unique extension to $\chi^{p,q}_c:=\sum_i
(-1)^i h^{p,q} H^i_c$, the $(p,q)$-th Euler functor of Deligne's mixed Hodge
structures for compactly supported cohomology \cite{Deligne2}, such that it
factors through the ring $K_0({\rm Var}_{\mathbb{C}})$:
$$
\xymatrix{{\rm Var}_{\mathbb{C}}\ar[rr]^{\chi_c^{p,q}}\ar[d]& &K_0({\rm Hodge})\cr
        K_0({\rm Var}_{\mathbb{C}})\ar@{-->}[urr]}
$$
Let $\mathbb{L}:=[\mathbb{A}^1_{\mathbb{C}}]$ be the Lefschetz class.
Notice that $\chi^{p,q}_c$ extends to $ K_0({\rm Var}_{\mathbb{C}})[\mathbb{L}^{-1}]$
since $\mathbb{L}$ corresponds to degree shifting operator of Hodge structures.\par

Now assume that $X$ is smooth and let $\phi: Y\to X$ be the blowing-up of $X$
along a smooth center $Z\subset X$ of codimension $r$, with exceptional divisor
$E\subset Y$. Then we have the well-known motivic equation for projective
bundles $E=\mathbb{P}_Z(N_{Z/X})\to Z$:
$$
[E] = [Z](1 + \mathbb{L} + \cdots \mathbb{L}^{r-1}) = [Z][\mathbb{P}^{r-1}].
$$
Since $[X] - [Z] = [Y] - [E]$, by formally localizing at $[\mathbb{P}^{r-1}]$
one gets
$$
[X] = [Y] - [E] + [Z] = ([Y] - [E]) + [E][\mathbb{P}^{r-1}]^{-1}.
$$
We call this a {\it nice change of variable formula} since the {\it Jacobian
factor} here depends only on the class $[E]$ instead of the precise structure
of the normal bundle $N_{Z/X}$. \par

The above computation can be performed inductively to show that, for $\phi:Y\to
X$ a composite of blowing-ups along smooth centers with $K_Y = \phi^* K_X +
\sum\nolimits_{i=1}^n e_i E_i$ and $E:=\bigcup_i E_i$ a normal crossing
divisor, the change of variable from $X$ to $Y$ reads
$$
[X] = \sum\nolimits_{I\subset \{1,\ldots,n\}}[E^\circ_I] \prod\nolimits_{i\in
I}[\mathbb{P}^{e_i +1}]^{-1},
$$
where $[E^\circ_I] := \bigcap_{i\in I}E_i \backslash \bigcup_{j\not\in I} E_j$.\par

According to the Meta Theorem, if we can prove such a change of variable
formula for any birational morphism $\phi: Y\to X$, we will then be able to
deduce that two $K$-equivalent smooth varieties $X$, $X'$ have $[X] = [X']$ in
$S^{-1}K_0({\rm Var}_{\mathbb{C}})$, where $S$ is the multiplicative set
generated by the classes of projective spaces. \footnote{I am grateful to W.\
Veys for suggesting the formulation of localizations on projective spaces.} Or
equivalently, $[P][X]=[P][X']$ for $P$ a product of projective spaces. Since
$\chi_c(V):=\sum_{p,q}\chi^{p,q}_c(V)$ is not a zero divisor for smooth
projective $V$, by applying the functor $\chi_c$ we conclude that $X$ and $X'$
have (non-canonically) isomorphic $\mathbb{Q}$-Hodge structures.
\par

The first proof of a weak form of this in certain completion of $K_0({\rm
Var}_{\mathbb{C}})[\mathbb{L}^{-1}]$ is due to Kontsevich and Denef-Loeser by
constructing motivic integration (see \S3.4 below). A numerical form of this
was also proved by Batyrev \cite{Batyrev2} using his version of motivic
integration. A full proof now is a consequence of the {\it Weak Factorization
Theorem} of Wlodarsczyk \cite{Wlodarsczyk} and Abramovich, Karu, Matsuki and
Wlodarsczyk \cite{AKMW}. It states that any birational map can be factorized
into sequences of blowing-ups and blowing-downs along smooth centers.\par

\subsection{Nash's Arc Spaces and the Motivic Integration}
During the same time the $p$-adic proof appeared, based on Kontsevich's idea,
Denef and Loeser \cite{DeLo1} had constructed the motivic integration over
Nash's formal arc spaces \cite{Nash}. In a correspondence with Loeser (c.f.\
\cite{Wang2}), we realized that instead of using $p$-adic integral, if we use
motivic integral then their {\it change of variable formula} (together with
Deligne's theory of mixed Hodge structures on arbitrary complex varieties as
above) will also lead to a proof of the above formula for arbitrary proper
birational morphism $\phi$ hence a proof of the strengthened result that the
Hodge numbers are the same. (See also \cite{Batyrev2} for an alternative
version of motivic integration.)\par

We give a brief sketch of their construction here. It starts with a measure
theory on the semi-algebraic subsets in Nash's arc spaces:
$$
\mu_X:\mathbb{B}(\mathcal{L}(X)) \longrightarrow
\widehat{K_0({\rm Var}_{\mathbb{C}})[\mathbb{L}^{-1}]}.
$$
Here $\mathcal{L}(X)$ as a set is simply the power series points of $X$ which
also has the structure as a pro-variety ${\rm ind.lim}\,\mathcal{L}_m(X)$ with
$\mathcal{L}_m(X)={\rm Hom}\,({\rm Spec}\,\mathbb{C}[t]/(t^{m+1}),X)$. There
are obvious natural morphisms $\pi_m:\mathcal{L}(X)\to\mathcal{L}_m(X)$ and
$\theta_m:\mathcal{L}_{m+1}(X)\to\mathcal{L}_m(X)$. When $X$ is smooth of
dimension $n$, $\pi_m$ is surjective and $\theta_m$ defines a piecewise affine
$\mathbb{A}^n_{\mathbb{C}}$ bundle structure on $\mathcal{L}_m(X)$. In general,
when $X$ is singular, $\pi_m$, $\theta_m$ are not even surjective. But we will
restrict to the smooth case here. \par

A set $S\subset \mathcal{L}(X)$ is stable (or cylindrical) if it is
semi-algebraic and is of the form $S = \pi_m^{-1}(A)$ for some (necessarily
constructible) set $A\subset \mathcal{L}_m(X)$. In this case we define its
$K_0({\rm Var}_{\mathbb{C}})[\mathbb{L}^{-1}]$-valued measure by $\mu_X(S) =
[\pi_m(A)]\mathbb{L}^{-mn}$. This has an unique extension to all semi-algebraic
sets $\mathbb{B}(\mathcal{L}(X))$ by partitioning them into certain countable
union of stable sets modulo some {\it measure zero sets}. The measure then
takes value in the completion of $K_0({\rm Var}_{\mathbb{C}})[\mathbb{L}^{-1}]$
with respect to the filtration $F^p := \{[S]\mathbb{L}^{-i}\,|\,\dim S - i \le
-p\}$. In particular $\mu_X(\mathcal{L}(X)) = [X]$ for smooth $X$. (Notice that
the motivic measure defined here differs from the one in \cite{DeLo1} by a
factor $\mathbb{L}^{-n}$.)
\par

Given a semi-algebraic set $S$ and a simple function $f:S\to
\mathbb{Z}\cup\{\infty\}$ such that $f^{-1}(k)$ is semi-algebraic all $k$, the
motivic integration is defined by
$$
\int_S \mathbb{L}^{-f} d\mu_X = \sum_{k\in \mathbb{Z}}\mathbb{L}^{-k}
\mu_X(f^{-1}(k)).
$$

Now for a birational morphism $\phi: Y\to X$ with $Y$ smooth, $\phi$ naturally
induces a map $\phi_*:\mathcal{L}(Y)\to\mathcal{L}(X)$. If $K_Y = \phi^*K_X +
E$ with $E$ a normal crossing divisor, the change of variable formula of Denef
and Loeser states that
$$
\int_{S}\mathbb{L}^{-f} d\mu_X
    = \int_{\phi^{-1}(S)}\mathbb{L}^{-f\circ\phi_*-{\rm ord}_t J\phi} d\mu_Y.
$$
Here $J\phi := \mathscr{O}_Y(-E)$ is the ideal sheaf generated by the
holomorphic Jacobian factor, ${\rm ord_t}\,\mathscr{I}: \mathcal{L}(X)\to
\mathbb{N}\cup\{0\}$ for any ideal sheaf $\mathscr{I}$ is the function of
minimal degree in $t$. Namely for $\gamma\in \mathcal{L}(X)$, ${\rm
ord_t}\,\mathscr{I}(\gamma) := {\rm min}_{g\in \mathscr{I}} \deg_t
g\circ\gamma(t)$. \par

The proof for general $X$ is technical. However for smooth $X$ the main idea of
the proof is not hard to explain. Indeed it is an application of the inverse
function theorem over power series rings which traces carefully the orders in
$t$. Since we only need the smooth case here, we will give an outline of the
proof in this case. Let $\phi:Y\to X$ be the birational morphism with $E_{\rm
red}\subset Y$ and $Z\subset X$ be the exceptional loci in $Y$ and $X$
respectively.\par

For each $k\in \mathbb{N}\cup\{0\}$ let $S_k\subset \mathcal{L}(Y)$ be the
subset $\gamma\in\mathcal{L}(Y)$ such that ${\rm ord}_t\,J\phi(\gamma) = k$. By
the inverse function theorem (see below), the map $\phi_*: \mathcal{L}(Y) \to
\mathcal{L}(X)$ is a bijection between $\mathcal{L}(Y)^\times : =
\mathcal{L}(Y)\backslash \mathcal{L}(E_{\rm red})$ and $\mathcal{L}(X)^\times:
= \mathcal{L}(X)\backslash \mathcal{L}(Z)$, thus there is no interesting
geometry on the map $\phi_*|_{S_k}:S_k \to \mathcal{L}(X)^\times$. However, the
important observation by Denef and Loeser is that when one takes finite
truncations in
$$
\xymatrix{ \mathcal{L}(Y)\ar[d]_\phi \ar[r]^{\pi_m} &
\mathcal{L}_m(Y)\ar[d]^{\phi_m}\\ \mathcal{L}(X)\ar[r]_{\pi_m} &
\mathcal{L}_m(X)}
$$
for all large enough $m$ the induced map $\phi_m|_{\pi_m(S_k)}:\pi_m(S_k) \to
\mathcal{L}_m(X)$ is indeed a piece-wise trivial $\mathbb{C}^k$ fibration over
its image. Together with the fact that $\mathcal{L}(Z)$ is {\it measure zero}
in $\mathcal{L}(X)$, this will imply the change of variable formula.\par

To investigate the fibration structure near one arc $\gamma\in \mathcal{L}(Y)$,
it is enough to restrict the map to formal neighborhoods
$\phi:\hat{\mathbb{C}}^n_{(0)}\to\hat{\mathbb{C}}^n_{(0)}$. Or equivalently to
represent $\phi$ by an algebraic map (still called $\phi$) on power series
$\phi:\mathbb{C}[[t]]^n\to\mathbb{C}[[t]]^n$ with $\phi(0)=0$. Let
$\phi(y(t))=x(t)$ with $y(t)\in S_k$ and let $\ell\ge 2k + 1$. We first notice
that for each $v\in\mathbb{C}[[t]]^n$, there is a unique solution
$u\in\mathbb{C}[[t]]^n$ of the equation
$$
\phi(y(t) + t^{\ell - k}u) = x(t) + t^\ell v.
$$
Indeed by Taylor's expansion
$$
\phi(y(t) + t^{\ell - k}u) = \phi(y(t)) + D\phi(y(t))t^{\ell-k}u + t^{2(\ell -
k)}R(t,u).
$$
Let $A = D\phi(y(t))$. The equation becomes $Au + R(t,u)t^{\ell-k} = t^kv$.
That is,
$$
u = (\det A)^{-1}t^kA^*(v - R(t,u)t^{\ell-2k}).
$$
Here $A^*$ is the adjoint matrix of $A$. Since ${\rm ord}_t\det A= {\rm ord}_t
J\phi(y(t)) = k$, the term $(\det A)^{-1}t^k$ has order zero. Also since $\ell
- 2k \ge 1$, by repeated substitutions this relation solves $u$ as a vector in
formal power series.\par

Now let $m\ge 2k$ and let $\ell = m + 1$. The above discussion shows that in
order to find all solutions of $\phi(\tilde{y}(t) \mod t^{m+1}) = x(t) \mod
t^{m+1}$, we may assume that $\tilde{y}(t) = y(t) + t^{m+1-k}u$. Notice that
the residue classes $\bar{u}=u \mod t^{k}$ form a linear space isomorphic to
$\mathbb{C}^{nk}$. By Hensel's lemma, in order to count the solutions we may
simply consider the equation $At^{m+1-k}\bar{u} = 0 \mod t^{m+1}$. That is,
$A\bar{u} = 0 \mod t^k$. Since ${\rm ord}_t\det(A^*)=(n-1)k$, the solution
space of $\bar{u}$ has dimension $nk - (n-1)k = k$ as expected.\par

This verifies that $\phi_m^{-1}\bar{x}(t)\cong\mathbb{C}^k$. The piece-wise
triviality needs other tools to prove it, which will not be reported here. For
the complete details the readers are referred to the original paper
\cite{DeLo1}.\par\smallskip

We remark that for $S=\mathcal{L}(X)$ and $E=\sum_{i=1}^n e_i E_i$ a normal
crossing, the change of variable formula gives $$ [X] =
\int_{\mathcal{L}(X)}\mathbb{L}^0 d\mu_X = \sum\nolimits_{I\subset
\{1,\ldots,n\}}[E^\circ_I] \prod\nolimits_{i\in I}
\frac{\mathbb{L}-1}{\mathbb{L}^{e_i +1}-1}. $$ Since
$\mathbb{L}^{e+1}-1=(\mathbb{L}-1)[\mathbb{P}^e]$, this coincides with the
formula in \S3.3.
\par

\section{$K$-equivalence Relation and Complex Elliptic Genera: Weak
Deformation/Decomposition Theorem}

\subsection{Some Background}
There is a build-in problem in all integration-theoretic approaches to the
$K$-equivalence relation. Namely we arrive at only $K$-theoretic or in practice
simply numerical conclusions. It is usually hard to get results of geometric
nature just from numerical data. In dimension three, the result of Koll\'ar and
Mori \cite{Kollar} on the flop decomposition of birational minimal models can
be easily generalized to any two $K$-equivalent threefolds. So the results
mentioned in \S1: naturally isomorphic cohomology groups, equivalent Hodge
structures and local moduli spaces are all still true. Moreover, these
canonical isomorphisms are all induced from the graph closure cycle of the
given birational map. It is clear that we can not achieve these statements from
integration theory only.\par

In the higher dimensional cases, due to the fact that it is (at least
currently) impossible to classify (terminal) singularities, the existence
problem of flops seems to be completely out of reach. This suggests that we
should not restrict the study of $K$-equivalence relation inside the category
of algebraic geometry only. We should allow (locally holomorphic) symplectic
deformations. That is, small deformations of almost complex structures which
are (integrable in a neighborhood of the exceptional loci and are) tamed by the
original symplectic form. In dimension three, with the help of classification
theory of singularities we may show that: if allowing symplectic deformations,
then any birational map between three dimensional $K$-equivalent manifolds can
be decomposed into composite of {\it classical flops} (see 4.2 below). All the
natural isomorphisms that we are interested in are then just simple
corollaries. The unsatisfactory fact is that we DO NOT know how to prove this
deformation/decomposition theorem directly without using the classification
theory. Such a proof should shed important light toward the higher dimensional
cases.\par

In fact, this symplectic deformation/decomposition theorem is even more useful
than the original flop decomposition theorem for certain problems. For example,
Li and Ruan \cite{LiRu} had shown in 1998 that it can be used to prove that
birational smooth minimal threefolds have equivalent quantum cohomology rings.
Notice that the ring structure of ordinary cohomology groups are not preserved
under flops $X\dashrightarrow X'$, in general $X$ and $X'$ are not even
homotopically equivalent.\par

\subsection{Some Well Known Flops}
For the reader's convenience, we recall the definition of certain known flops.
The simplest type of flops are called {\it ordinary flops}. An ordinary
$\mathbb{P}^r$-flop (or simply $\mathbb{P}^r$-flop) $f:X\dashrightarrow X'$ is
a birational map such that the exceptional set $Z\subset X$ has a
$\mathbb{P}^r$-bundle structure $\psi:Z\to S$ over some smooth variety $S$ and
the normal bundle $N_{Z/X}$ is isomorphic to $\mathscr{O}(-1)^{r+1}$ when
restricting to any fiber of $\psi$. The map $f$ and the space $X'$ are then
obtained by first blowing up $X$ along $Z$ to get $Y$, with exceptional divisor
$E$ a $\mathbb{P}^r\times\mathbb{P}^r$-bundle over $S$, then blowing down $E$
along another fiber direction. Ordinary $\mathbb{P}^1$-flops are also called
classical flops. Three dimensional classical flops are the most well-known
Atiyah flops over $(-1,-1)$ rational curves. \par

Another important example is the {\it Mukai flops} $f:X\dashrightarrow X'$. In
this case it is required that the exceptional set $Z\subset X$ is of
codimension $r$ and has a $\mathbb{P}^r$-bundle structure
$\psi:Z=\mathbb{P}_S(F)\to S$ (for some rank $r+1$ vector bundle $F$) over a
smooth base $S$, moreover the normal bundle $N_{Z/X}\cong T^*_{Z/S}$, the
relative cotangent bundle of $\psi$. To get $f$, one first blows up $X$ along
$Z$ to get $\phi:Y\to X$ with exceptional divisor
$E=\mathbb{P}_Z(T^*_{Z/S})\subset \mathbb{P}_S(F) \times_S\mathbb{P}_S(F^*)$ as
the incidence variety. The first projection corresponds to $\phi$ and one may
contract $E$ through the second projection to get $\phi':Y\to X'$. Mukai flops
naturally occur in hyperk\"ahler manifolds \cite{Huyb1}.\par

\subsection{Main Conjectures}
In 2000, the author made a series of conjectures on $K$-equivalent manifolds
$X$ and $X'$ under birational map $f:X \dashrightarrow X'$:
\begin{itemize}
\item[I.] The morphism $T:H^k(X,\mathbb{Q})\to H^k(X',\mathbb{Q})$
induced from the graph closure $\bar\Gamma_f\subset X\times X'$ is an
isomorphism which preserves the rational Hodge structures. There also exists a
canonical correspondence $\bar\Gamma_f + \sum_i T_i\subset A^n(X\times X')$,
with $T_i$ being certain degenerate correspondences, which defines an
isomorphism on integral cohomology groups modulo torsion. \footnote{I am
grateful to D.\ Huybrechts and Y.\ Namikawa for pointing out the necessity to
modify the graph closure in order to get the conjectural isomorphisms on
integral cohomologies.}
\item[II.] The local deformation spaces ${\rm Def}(X)$ and ${\rm Def}(X')$ are
canonically isomorphic in the sense that the local universal families are
$K$-equivalent over the base. Moreover, suitable compactifications of their
polarized moduli spaces should also be $K$-equivalent.
\item[III.] $X$ and $X'$ have canonically isomorphic quantum
cohomology rings over the extended K\"ahler moduli spaces. In other words,
their quantum cohomology rings can be analytically continued to each other.
\item[IV.] Deformation/Decomposition Theorem: under generic symplectic
perturbations which respect $f$, the deformed $f$ can be decomposed into finite
copies of ordinary $\mathbb{P}^r$-flops for various $r$'s.
\end{itemize}
It is also expected that IV would be the key step toward resolving conjectures
I, II and III.\par

The main progress made in \cite{Wang3} is to prove a weak form of conjecture
IV: if we further modulo complex cobordism, then any birational map between
$K$-equivalent manifolds can be decomposed into the composite of finite number
of ordinary $\mathbb{P}^1$-flops. Notice that since in a flat family of
algebraic cycles the dimension can not go down under specialization, the $r$'s
appear in conjecture IV can not take the value $1$ only. This explains that the
weak form we proved is still far away from the original conjecture. Another
important remark is related to the Mukai flops. These flops are not included in
conjecture IV since Huybrechts \cite{Huyb1} had shown in 1996 that Mukai flops
in hyperk\"ahler manifolds will disappear (become isomorphisms) under generic
deformations. Recently he completed the discussion by showing that birational
hyperk\"ahler manifolds become isomorphic under generic deformations of complex
structures \cite{Huyb2}. Huybrechts' results can be regarded as one of the most
important evidences of the above conjectures (c.f.\ 5.3).\par

Also it should be remarked that for Calabi-Yau manifolds, the equivalence of
Hodge numbers gives numerical evidence for Main Conjecture II since the
relevant groups in the Kodaira-Spencer theory are all Hodge groups:
$$
H^i(X, T_X) \cong H^i(X, \Omega^{n-1}_X) \cong H^{n-1,i}(X).
$$
However, in order to proceed, we really need the validity of Conjecture I.\par

\subsection{Complex Elliptic Genera under $\mathbb{P}^1$-flops}
As for the proof of the weak deformation/decomposition theorem, we notice that
according to a result of Milnor \cite{Milnor} and Novikov \cite{Novikov}, the
complex cobordism class of a compact stably almost complex manifold ($X$ such
that $T_X\oplus \xi$ has a complex vector bundle structure for some trivial
bundle $\xi$) is determined precisely by all its Chern numbers. The complex
cobordism ring $\Omega^U$ is defined to be the ring of compact stably almost
complex manifolds modulo cobordism by such manifolds with boundaries.\par

Recently Totaro \cite{Totaro} showed that the most general Chern numbers that
are invariant under ordinary $\mathbb{P}^1$-flops consists of the so-called
{\it complex elliptic genera}. Recall that an $R$-genus is nothing but a ring
homomorphism $\varphi:\Omega^U\to R$. Equivalently it can be defined through
Hirzebruch's power series recipe \cite{Hirz1}. Let $Q(x)\in R[\![x]\!]$ and $X$
be an almost complex manifold with formal Chern roots decomposition
$c(T_X)=\prod_{i=1}^n (1+x_i)$, then
$$
\varphi_Q(X) :=\prod\nolimits_{i=1}^n Q(x_i)[X]=:\int_X K_Q(c(T_X))
$$
defines an $R$-genus. As usual write $Q(x)=x/f(x)$ then the complex elliptic
genera is defined by the three parameter ($k\in\mathbb{C}$, $\tau$ and a marked
point $z$) power series
$$
f(x)=e^{(k+\zeta(z))x}\,\frac{\sigma(x)\sigma(z)}{\sigma(x+z)}.
$$

Hirzebruch \cite{Hirz2} has reproved Totaro's theorem using Atiyah-Bott
localization theorem. He showed that $\varphi_Q$ is invariant under
$\mathbb{P}^1$-flops if and only if $F(x):=1/f(x)$ satisfies the functional
equation $$ F(x+y)(F(x)F(-x)-F(y)F(-y))=F'(x)F(y)-F'(y)F(x). $$ Moreover, the
solutions is given by the above $f$ exactly.\par

\subsection{Complex Elliptic Genera under $K$-equivalence}
The main contribution in \cite{Wang3} is to show that complex elliptic genera
are also invariant among general $K$-equivalent manifolds. Hence in $\Omega^U$,
the ideal $I_1$ generated by $[X]-[X']$ with $X$ and $X'$ related by a
$\mathbb{P}^1$-flop is equal to the seemingly much lager ideal $I_K$ generated
by $K$-equivalent pairs $[X]-[X']$.

Following the meta theorem, the most important step in the proof is to develop
a change of variable formula for genera (or Chern numbers) under blowing-ups.
First, using standard intersection theory and Hirzebruch's theory of virtual
genus \cite{HBJ}, we proved a residue formula for a single blowing-up
$\phi:Y\to X$ along smooth center $Z$ of codimension $r$. Namely, for any power
series $A(t)\in R[\![t]\!]$:
\begin{eqnarray*}
\int_Y A(E)\,K_Q(c(T_Y))
 \!\!&=&\!\! \int_X A(0)\,K_Q(c(T_X))\\
 && + \int_Z {\rm Res}_{\,t=0}\Big(\frac{A(t)}{f(t)\prod_{i=1}^r f(n_i - t)}\Big)
 \,K_Q(c(T_Z)).
\end{eqnarray*}
Here $n_i$'s denote the formal Chern roots of the normal bundle $N_{Z/X}$ and
the residue stands for the coefficient of the degree $-1$ term of a Laurent
power series with coefficients in the cohomology ring or the Chow ring of $X$.
In order to have a change of variable formula, we need the residue term to
vanish. If we already know the expression of $f$ as above, than for $z$ not an
$r$-torsion point it is not hard to find $$ A(t,r) =
e^{-(r-1)(k+\zeta(z))t}\,\frac{\sigma(t + rz)\sigma(z)}
{\sigma(t+z)\sigma(rz)}. $$ In fact the $r=2$ case corresponds to a functional
equation $$ \frac{1}{f(x)f(y)} = \frac{A(x)}{f(x)f(y - x)} +
\frac{A(y)}{f(y)f(x - y)} $$ which also has solutions given by $f$ (and $A$ is
determined by $f$), but with $z$ not a 2-torsion points. \footnote{It is
expected that complex genera which admit the change of variable formula for
codimension $r$ blowing-ups consists of precisely the complex elliptic genera
with $z$ not an $r$-torsion point.} Notice that a classical theorem of
Weierstrass states that solutions of functional equations involving only
$f(x)$, $f(y)$ and $f(x + y)$ are constructed from the Weierstrass elliptic
functions. However, the functional equations appeared here (and also the one
considered by Hirzebruch in \S4.4) are not of this type.

The general change of variable reads: let $\varphi=\varphi_{k,\tau,z}$ be the
complex elliptic genera. Then for any algebraic cycle $D$ in $X$ and birational
morphism $\phi:Y \to X$ with $K_Y = \phi^*K_X + \sum e_i E_i$, we have (write
$K_Q=K_{\varphi}$)
$$
\int_D K_\varphi(c(T_X))
= \int_{\phi^*D} \prod\nolimits_i A(E_i, e_i + 1)\,K_\varphi(c(T_Y)).
$$
Or equivalently, $\phi_*\prod_i A(E_i, e_i +
1)\,K_\varphi(c(T_Y))=K_\varphi(c(T_X))$. This is first proved by induction for
$\phi$ a composite of blowing-ups. The general cases can be reduced to the
blowing-up case by applying the weak factorization theorem \cite{Wlodarsczyk}
\cite{AKMW}. (A similar result for the case $k=0$ (elliptic genera) was also
obtained recently by Borisov and Libgober \cite{BoLi}.)\par

The formula implies that $X=_K X'\Rightarrow \varphi_{k,\tau,z}(X)=
\varphi_{k,\tau,z}(X')$ with $z$ not a torsion point. But then we also get
$\varphi(X)=\varphi(X')$ in all cases by continuity.\par

Notice that it is symbolically convenient to denote $K_\varphi(c(T_X))$ by
$d\mu_X$ and regard it as an {\it elliptic measure}, though we do not really
construct a measure theory as in \S3. When we specialize to Todd genus
(rational measure), the Jacobian factor reduces to $A = 1$ and the change of
variable formula is a simple corollary of the Grothendieck Riemann-Roch
theorem.
\par

It is worth mentioning that except for the last step, the proof works both in
the category of (stably almost) complex manifolds and in the category of
algebraic manifolds in arbitrary characteristic. While the weak factorization
theorem has been proved for both the complex analytic category \footnote{I am
grateful to B.\ Totaro for pointing out this to me.} and also the algebraic
category, the later is proved under the restriction over fields of
characteristic zero. It is expected that one should find a Grothendieck
Riemann-Roch type argument to replace the weak factorization theorem to get a
more satisfactory conclusion.\par

\subsection{Chern Numbers versus Hodge Numbers}
Hodge numbers and Hodge structures determine a substantial part of the complex
elliptic genera and also give information to the complex moduli. Recall that
\cite{Totaro}:
$$
\varphi(X) = \chi\Big(X, K_X^{\otimes(-k)}\otimes\prod\nolimits_{m\ge 1}
(\Lambda_{-y^{-1}q^m}T\otimes\Lambda_{-y^{-1}q^{m-1}}T^*\otimes S_{q^m}T\otimes
S_{q^m}T^*)\Big)
$$
for $q = e^{2\pi i\tau}$, $y=e^{2\pi iz}$ and $T = T_X -n$ the rank zero virtual
tangent bundle. The twisted $\chi_y$-genus corresponds to the two parameter genera
$$
\chi_y(X) := \chi\Big(X,K_X^{\otimes(-k)}\otimes\Lambda_y T_X^*\Big),
$$
which is equivalent to knowing all $\chi(X,K_X^{\otimes(-k)}\otimes\Omega_X^p)$
for $p\ge 0$. If $n=\dim X\le 11$, the twisted $\chi_y$ genus contains the same
Chern numbers as the complex elliptic genera. So in this range, twisted
$\chi_y$ genus contains precisely all Chern numbers that are invariant under
the $K$-equivalence relation. If $n\le 4$, the twisted $\chi_y$ genus contains
all Chern numbers, so all Chern numbers are invariant under $K$-equivalence for
dimensions up to $4$.\par

It is clear that if $K_X$ is trivial, that is, $X$ is a Calabi-Yau manifold,
then the twisted $\chi_y$ genus becomes Hirzebruch's $\chi_y$ genus $\sum_{p\ge
0}\chi(X,\Omega_X^p)\,y^p$. In particular, it is determined by the Hodge
numbers. So the equivalence of elliptic genera (that is, $k=0$) follows from
the equivalence of Hodge numbers when $n\le 11$. But when $n\ge 12$, the
elliptic genera and Hodge numbers contain quite a different type of
information.\par

\section{Other Aspects of $K$-equivalence Relation}

For completeness, we add a few topics that are closely related to the study of
$K$-equivalence relation but not directly related to the author's current
approaches. The interested reader should consult the original papers for more
details.\par

\subsection{$K$-equivalence for Singular Varieties}
The notion of $K$-partial ordering makes sense for general
$\mathbb{Q}$-Gorenstein varieties. For log-terminal varieties, the integration
formalism works well as in the smooth case (the measure is finite
$\Leftrightarrow$ $X$ is log-terminal) and $K$-equivalence still implies
measure-equivalence, both the $p$-adic and motivic ones. The major problem here
is to understand the geometric meaning of the total measure. For $p$-adic
measure, it is a weighted counting of rational points over finite fields, but
we do not know how to make it precise. \par

In \cite{Batyrev2} Batyrev defined the stringy $E$ function for log-terminal
varieties and also the {\it stringy Hodge numbers} when this $E$ function is a
polynomial. In terms of motivic measure, it can be defined by taking
$h^{p.q}_{\rm st}(X)=\chi^{p,q}(\mu_X(\mathcal{L}(X))$. It is clear that
$K$-equivalent varieties have the same stringy Hodge numbers, but its meaning
still needs to be further clarified. Does there exist corresponding cohomology
theories (spaces)? Noticed that $h^{p,q}_{\rm st}(X)$ is in general only a
rational number and may as well be negative. Veys has recently investigated the
situation for normal surface singularities \cite{Veys}.\par

A more manageable case is the crepant resolutions of Gorenstein quotient
singularities $\phi:Y\to X:=\mathbb{C}^n/G$ with $G\subset {\rm
SL}(n,\mathbb{C})$ a finite subgroup such that $K_Y=\phi^* K_X$. The {\it McKay
correspondence} asserts that, among other things, a natural basis of $H^*(Y)$
is in one to one correspondence with conjugacy classes of $G$. The numerical
version of it has been proved by Batyrev \cite{Batyrev3} and also by Denef and
Loeser \cite{DeLo2} using motivic integration. The geometric correspondence has
not been treated except in dimensions $\le 3$ \cite{BKR}. Moreover, it is not
known in general when does $X$ admit crepant resolutions. A possible approach
is to look at the $p$-adic measure of $X$. As it must be the measure of $Y$
when $Y$ exists, the corresponding counting function must behave like counting
points on smooth varieties. This should put certain restrictions on $G$. \par

A recent attempt toward constructing the expected cohomology theory was given
by Chen and Ruan's {\it orbifold cohomology} for Gorenstein orbifolds
\cite{ChRu}. But the naturality problem has not been solved yet. It seems to be
a difficult problem for the construction for general $\mathbb{Q}$-Gorenstein
varieties. The author do not know even conceptually how to extend the Main
Conjectures to the singular case.\par

\subsection{Derived Categories and Fourier-Mukai Transform}

In 1995, Bondal and Orlov \cite{BoOr} showed that for special cases of
$\mathbb{P}^r$-flops (namely $S$ is a point in the notation of 4.2), The
natural transform $\Phi := {\bf R}\phi'_*\circ{\bf L}\phi^*$ gives an
equivalence of triangulated categories $D^b(X)\cong D^b(X')$ (here $D^b(X)$
denotes the bounded derived category of coherent sheaves). Later on Bridgeland
\cite{Brid} had made significant progress on this approach by extending their
result to all smooth threefold flops. The important issue here is that the
flopped variety $X^+$ may be constructed as certain {\it fine} moduli spaces.
More precisely, he showed that, for $\psi:X\to\bar X$ a flopping contraction
(that is, $K_X$ is $\psi$-trivial) from a smooth threefold $X$, let $X^+:={\rm
Per}(X/\bar X)$ be the distinguished component of the moduli space of {\it
perverse point sheaves}, then $X^+$ is smooth and $f:X\dashrightarrow X^+$ is
the flop. Also the Fourier-Mukai transform
$$
{\bf R}{p_2}_*(\mathcal{E}\mathop{\otimes}\limits^{L}{\bf L}p_1^*(-)):
D^b(X)\to D^b(X^+)
$$
(here $p_1$, $p_2$ are the projections from $X\times X^+$ to $X$ and $X^+$
respectively) induced by the {\it universal perverse point sheave}
$\mathcal{E}\in D^b(X\times X^+)$ is an equivalence of categories.\par

Bridgeland's theorem is recently generalized by Chen \cite{Chen} to 3-folds
with Gorenstein terminal singularities and by Kawamata \cite{Kawa3} to three
dimensional orbifolds. By combining with Chen's result, Kawamata also proved
the equivalence of derived categories for all three dimensional terminal flops
\cite{Kawa4}. There seems to be of some hope to deal with certain higher
dimensional flopping contractions $\psi:X\to \bar X$ with relative dimension
$\le 1$ through their methods. \par

In \cite{Kawa4}, Kawamata conjectured that for birational projective manifolds,
the notion of $K$-equivalence should be equivalent to $D$-equivalence, namely
varieties with equivalent derived categories of coherent sheaves. This is
clearly closely related to our main conjectures, but a precise relation between
derived categories and cohomologies does not seem to be well studied yet.

\subsection{Flop Decomposition for Hyperk\"ahler Manifolds}
Hyperk\"ahler manifolds (or holomorphically symplectic manifolds) have been
extensively studied lately. All our Main Conjectures follow from Huybrechts'
fundamental works \cite{Huyb1} \cite{Huyb2} mentioned in \S4.3. For Conjecture
I the correspondence cycle $\Gamma\subset X\times X'$ used by him is the
limiting cycle $\lim_{t\to 0}\bar\Gamma_{f_t}$ induced from nearby isomorphisms
$f_t:\mathscr{X}_t\cong\mathscr{X}'_t$ with $t\ne 0$. This cycle in general
contains more than one irreducible components. In fact for a Mukai flop, the
map $T$ induced from the graph closure will in general preserve only the {\it
rational} cohomologies. The statement in Conjecture I is still unknown for
birational hyperk\"ahler manifolds under the map $T$.
\par

On the other direction, Burns, Hu and Luo \cite{BHL} had shown that birational
maps between hyperk\"ahler fourfolds can be decomposed into composite of Mukai
flops, if all the irreducible components of the exceptional loci are normal.
Very recently this normality assumption was justified by Wierzba and
Wi\'sniewski \cite{WW}, hence the four dimensional case was settled completely.
Notice that in this case the Mukai flop is of a particularly simple type (in
the notation of \S4.2, $Z\cong\mathbb{P}^2$ and $S$ is a single point). Since
the equivalence of derived categories for Mukai flops with $S$ being a point is
proved by Namikawa \cite{Namikawa} and Kawamata \cite{Kawa4}, we see that
birational hyperk\"ahler fourfolds are indeed $D$-equivalent. \par

As a final question, can one prove the above results on $D$-equivalence without
making use of the explicit flops decomposition? Notice that this is the main
theme of our approach toward cohomologies and complex genera in higher
dimensions.

\bibliographystyle{amsplain}

\begin{thebibliography}{10}

\bibitem{AKMW} D. Abramovich, K.Karu, K. Matsuki and J. Wlodarsczyk;
{\it Torifications and factorizations of birational maps}, J.\ Amer.\ Math.\
Soc.\ {\bf 15} (2002), 531-572.
\bibitem{Batyrev1} V. Batyrev; {\it On the Betti numbers of birationally
isomorphic projective varieties with trivial canonical bundles}, {\tt
alg-geom}/9710020.
\bibitem{Batyrev2} V. Batyrev; {\it Stringy Hodge numbers of varieties with
Gorenstein canonical singularities}, in {\it Integrable systems and algebraic
geometry} (Kobe/Kyoto, 1997), 1-32, World Sci.\ Publishing, Rivers Edge, NJ,
1998.
\bibitem{Batyrev3} V. Batyrev; {\it Non-Archimedian integral and stringy Euler
numbers of log terminal pairs}, J.\ Eur.\ Math.\ Soc.\ {\bf 1} (1999), no.\ 1,
5-33.
\bibitem{BPV} W. Barth, C. Peters and A. Van de Ven; {\it Compact Complex
Surfaces}, Springer Verlag 1984.
\bibitem{BoOr} A. Bondal and D. Orlov; {\it Semiorthogonal decompositions of
algebraic varieties}, {\tt math.AG}/9506012.
\bibitem{BoLi} L. Borisov and A. Libgober; {\it Elliptic genera of singular
varieties}, {\tt math.AG}/0007108.
\bibitem{Brid} T. Bridgeland; {\it Flops and derived categories};
{\tt math.AG}/0009053.
\bibitem{BKR} T. Bridgeland, A, King and M. Reid; {\it Mukai implies McKay};
J.\ Amer.\ Math.\ Soc.\ {\bf 14} (2001), 535-554.
\bibitem{BHL} D. Burns, Y. Hu and T. Luo; {\it Hyperk\"ahler manifolds and
birational transformations in dimension 4}, {\tt math.AG}/0004154.
\bibitem{Chen} J.-C. Chen; {\it Flops and equivalence of derived categories
for threefolds with only terminal Gorenstein singularities}, {\tt math.AG}/0202005.
\bibitem{ChRu} W. Chen and Y. Ruan; {\it A new cohomology theory for orbifold},
{\tt math.AG}/0004129.
\bibitem{Deligne1} P. Deligne, {\it La conjecture de Weil I}, IHES Publ.\ Math.\
{\bf 43} (1974), 273-307. {\it II}, {\bf 52} (1980), 313-428.
\bibitem{Deligne2} P. Deligne. {\it Th\'eorie de Hodge II}, IHES Publ.\ Math.\
{\bf 40} (1971), 5-57. {\it III}, {\bf 44} (1974), 5-77.
\bibitem{DeLo1} J. Denef and F. Loeser; {\it Germs of arcs on singular algebraic
varieties and motivic integration}, Inv.\ Math.\ {\bf 135} (1999), 201-232.
\bibitem{DeLo2} J. Denef and F. Loeser; {\it Motivic integration, quotient
singularities and the McKay correspondence}, {\tt math.AG}/9903067.
\bibitem{Faltings} G. Faltings; {\it $p$-adic Hodge thoery},
J.\ Amer.\ Math.\ Soc.\ {\bf 1} No.1 (1988), 255-299.
\bibitem{FoMe} J.-M. Fontaine and W. Messing; {\it $p$-adic periods and
$p$-adic \'etale cohomology}, Contemp.\ Math.\ {\bf 67} (1987), 179-207.
\bibitem{Freidman} R. Freidman, {\it Simultaneous resolution of threefold double
points}, Math.\ Ann.\ {\bf 274} (1986) 671--689.
\bibitem{Fulton} W. Fulton; {\it Intersection Theory}, Erge.\ Math.\ ihr.\ Gren.;3.
Folge, Bd 2, Springer-Verlag 1984.
\bibitem{Hirz1} F. Hirzebruch; {\it Topological Methods in Algebraic Geometry},
Grund.\ der Math.\ Wissen.\ {\bf 131}, Springer 1966.
\bibitem{Hirz2} F. Hirzebruch; {\it Complex cobordism and the elliptic genus},
Contemp.\ Math.\ {\bf 241} (1999), 9-20.
\bibitem{HBJ} F. Hirzebruch, T. Berger and R. Jung; {\it Manifolds and Modular Forms},
Max-Plank-Institut f\"ur Math.\ 1992.
\bibitem{Huyb1} D. Huybrechts; {\it Birational symplectic manifolds and their
deformations}, J.\ Diff.\ Geom.\ {\bf 45} (1997), 488-513.
\bibitem{Huyb2} D. Huybrechts; {\it Compact hyperk\"ahler manifolds: basic
results}, Invent.\ Math.\ {\bf 135} (1999), 63-113. Erratum {\tt math.AG}/0106014.
\bibitem{Ito} T. Ito; {\it Stringy Hodge numbers and $p$-adic Hodge theory},
preprint UTMS-2002-1, Univ.\ of Tokyo.
\bibitem{Kawa1} Y. Kawamata; {\it The crepant blowing-up of 3-dimensional
canonical singularities and its applications to degenerations of surfaces};
Annals of Math.\ {\bf 127} (1988), 93-163.
\bibitem{Kawa2} Y. Kawamata; {\it Abundance theorem for algebraic threefolds},
Invent.\ Math.\ {\bf 108} (1992), 229-246.
\bibitem{Kawa3} Y. Kawamata; {\it Francia's flip and derived categories};
{\tt math.AG}/0111041.
\bibitem{Kawa4} Y. Kawamata; {$D$-equivalence and $K$-equivalence}, preprint
2002, {\tt math.AG}/0205287.
\bibitem{KMM} Y. Kawamata, K. Matsuda, K. Matsuki; {\it Introduction to
the minimal model program}, Adv.\ Stud.\ in Pure Math.\ {\bf 10} (1987), 283-360.
\bibitem{Kollar} J. Koll\'ar; {\it Flops}, Nagoya Math.\ J. {\bf 113} (1989), 15--36.
\bibitem{KoMo1} J. Koll\'ar and S. Mori; {\it Classification of three dimensional
flips}, J.\ Amer.\ Math.\ Soc.\ {\bf 5} No.3 (1992), 533-703.
\bibitem{KoMo2} J. Koll\'ar and S. Mori; {\it Birational Geometry of Algebraic
Varieties}, Cambridge University Press 1998.
\bibitem{LiRu} A.-M. Li and Y. Ruan; {\it Symplectic surgery and Gromov-Witten
invariants of Calabi-Yau 3-folds}, Invent.\ Math.\ {\bf 145} (2001), 151-218.
\bibitem{Mori1} S. Mori; {\it Threefolds whose canonial bundles are not
numerically effective}, Annals of Math.\ {\bf 116} (1982), 133-176.
\bibitem{Mori2} S. Mori; {\it On 3-dimensional terminal singularities},
Nagoya Math.\ J.\ {\bf 98} (1985), 43-66.
\bibitem{Mori3} S. Mori; {\it Flip theorem and minimal models for 3-folds},
J.\ Amer.\ Math.\ Soc.\ {\bf 1} (1988), 117-253.
\bibitem{Milnor} J. Milnor; {\it On the cobordism ring $\Omega^*$ and a complex
analogue}, Amer.\ J.\ Math.\ {\bf 82} (1960), 505-521.
\bibitem{Namikawa} Y. Namikawa; {\it Mukai flops and derived categories},
preprint 2002, {\tt math.AG}/0203287.
\bibitem{Nash} J.F. Nash; {\it Arc structures of singularities}, preprint 1968,
Duke Math.\ J.\ {\bf 81} (1995), 31-38.
\bibitem{Novikov} S. Novikov; {\it Homotopy properties of Thom complexes},
Mat.\ Sb.\ {\bf 99} (1962), 407-442.
\bibitem{Reid} M. Reid, {\it Young person's guide to canonical singularities},
in Algebraic Geometry Bowdowin 1985, Proc.\ Symp.\ Pure Math.\ {\bf 46} (1987).
\bibitem{Serre} J.-P. Serre; {\it Abelian $\ell$-adic Representations and
Elliptic Curves}, W.\ A.\ Benjamin, Inc.\ 1968.
\bibitem{Totaro} B. Totaro; {\it Chern numbers for singular varieties and
elliptic homology}, Annals of Math.\ {\bf 151} 2000, 757-791.
\bibitem{Veys} W. Veys; {\it Stringy invariants of normal surfaces}, {\tt
math.AG}/0205293.
\bibitem{Wang1} C.-L. Wang; {\it On the incompleteness of the Weil-Petersson metric
along degenerations of Calabi-Yau manifolds}, Math. Res. Let {\bf 4}
(1997), 157-171.
\bibitem{Wang2} C.-L. Wang; {\it On the topology of birational minimal models},
J.\ Diff.\ Geom.\ {\bf 50} (1998), 129-146.
\bibitem{Wang3} C.-L. Wang; {\it $K$-equivalence in birational geometry and
characterizations of complex elliptic genera}, to appear in J.\ Alg.\ Geom.\ 2002.
\bibitem{Wang4} C.-L. Wang; {\it Cohomology theory in birational geometry},
Report on NSC Project 89-2115-M-002-012, Taiwan 2000. To appear.
\bibitem{Weil} A. Weil; {\it Ad\`ele and Algebraic Groups}, Prog.\ Math.\
{\bf 23}, Birkhauser, Boston 1982.
\bibitem{WW} J. Wierzba and J.A. Wi\'sniewski; {\it Small contractions of
symplectic 4-folds}, {\tt math.AG}/0201028.
\bibitem{Wilson} P.M.H. Wilson; {\it Symplectic deformations of Calabi-Yau 3-folds},
J.\ Diff.\ Geom.\ {\bf 45} (1997), 611-637.
\bibitem{Wlodarsczyk} J. Wlodarsczyk; {\it Combinatorial structures on toroidal
varieties and a proof of the weak factorization theorems}, preprint 1999,
{\tt math.AG}/9904076.
\bibitem{Yau} S.T.\ Yau, {\it On the Ricci curvature of a compact K\"ahler
manifold and the complex Monge-Amp\`ere equation I}, Comm.\ Pure and
Appl.\ Math.\ {\bf 31} (1978), 339--441.
\end{thebibliography}

\end{document}